\documentclass[11pt]{article}
\usepackage{graphicx}
\usepackage{amsmath}
\usepackage{color}
\usepackage{xcolor}
\usepackage{authblk}
\usepackage[english]{babel}
\usepackage{amssymb,amsfonts,amsthm}
\usepackage{enumitem}

\newif\ifpdf
\ifx\pdfoutput\undefined \pdffalse \else \pdfoutput=1 \pdftrue \fi
\ifpdf \else \fi \textwidth = 5.5in \textheight = 8.5 in
\newtheorem{theorem}{Theorem}

\newtheorem{corollary}[theorem]{Corollary}

\newtheorem{definition}[theorem]{Definition}

\newtheorem{lemma}[theorem]{Lemma}

\newtheorem{observation}[theorem]{Observation}
\newtheorem{problem}[theorem]{Problem}
\newtheorem{proposition}[theorem]{Proposition}

\makeatletter
\def\imod#1{\allowbreak\mkern10mu({\operator@font mod}\,\,#1)}
\makeatother

\begin{document}

\title{Majority dominator colorings of graphs}

\author[1]{Marcin Anholcer}
\author[2]{Azam Sadat Emadi}
\author[2]{Doost Ali Mojdeh}

\affil[1]{\scriptsize{}Pozna\'n University of Economics and Business, Institute of Informatics and Quantitative Economy}
\affil[ ]{Al.Niepodleg{\l}o\'sci 10, 61-875 Pozna\'n, Poland, \textit{m.anholcer@ue.poznan.pl}}
\affil[ ]{}
\affil[2]{Department of Mathematics, Faculty of Mathematical Sciences,  University of Mazandaran}
\affil[ ]{Babolsar, Iran, \textit{math\_emadi2000@yahoo.com, damojdeh@umz.ac.ir}}

\maketitle

\begin{abstract}
Let $G$ be a simple graph of order $n$. A majority dominator coloring of a graph $G$ is proper coloring in which each vertex of the graph dominates at least half of one color class. The majority dominator chromatic number $\chi_{md}(G)$ is the minimum number of color classes in a majority dominator coloring of $G$. In this paper we study properties of the majority dominator coloring of a graph. We obtain tight upper and lower bounds in terms of chromatic number, dominator chromatic number, maximum degree, domination and independence number. We also study majority dominator coloring number of selected families of graphs.
\end{abstract}
\textbf{Keywords}: Majority dominator chromatic number, majority dominator coloring, chromatic number, independence number, domination number.\\
\textbf{MSC 2010:} 05C15, 05C78, 05C69


\section{Introduction}

Let $G = (V, E)$ be a simple graph. For any vertex $ v\in V$, the open neighborhood of $v$ is the set $N(v)=\{u\in V|uv\in E\} $ and the closed neighborhood of $v$ is the set $N[v] = N(v)\cup \{v\}$. For a set $S\subseteq V$, the open neighborhood of $S$ is $N (S)=\bigcup_{v\in S}N(v)$ and the closed neighborhood of $S$ is $N[S]=N(S)\cup S$. A set $S\subseteq V$ is a dominating set if $N[S] = V$, or equivalently, every vertex in $V\setminus S$ is adjacent to at least one vertex in $S$. The domination number $\gamma(G)$ is the minimum cardinality of a dominating set in $G$. A dominating set with cardinality $\gamma(G)$ is called a $\gamma(G)$-set, for more information of these parameters, refer to \cite{bl, bz}. A set $I\subseteq V(G)$ is called independent set if the induced subgraph $G[I]$ has no edge. The size of maximum independent set in $G$ is called independence number and is denoted by $\alpha(G)$ \cite{west}. Finally, the size of any maximum matching in $G$ is called the matching number and denoted with $\nu(G)$.

Further in this paper we are going to consider some graph products. Let $G$ and $H$ be two graphs. The \emph{strong product} $G\boxtimes H$ is a graph with vertex set $V(G)\times V(H)$. Two vertices $(g,h)$ and $(g^{\prime},h^{\prime })$ are adjacent in $G\boxtimes H$ if either $g=g^{\prime }$ and $h$ is adjacent with $h^{\prime }$ in $H$, or $h=h^{\prime }$ and $g$ is adjacent with $g^{\prime }$ in $G$, or $g$ is adjacent with $g^{\prime }$ in $G$ and $h$ is adjacent with $h^{\prime }$ in $H$ (see e.g. \cite[p.36]{IK}). The corona product $G\circ H$ is in turn a graph consisting of $G$ and $|V(G)|$ copies of $H$, where each copy $H_j$ corresponds with a different vertex $v_j\in V(G)$ and all its vertices are adjacent with $v_j$.

Proper coloring of a graph $G$ is an assignment of colors to its vertices, such that no two adjacent vertices receive the same color. The chromatic number $\chi(G)$ is the minimum number of colors required for a proper coloring of $G$. Dominator coloring is a proper coloring in which each vertex of the graph dominates (i.e., is adjacent to) every vertex of some color class (it can be its own color class, in particular when it consists of only one vertex). The dominator chromatic number $\chi_{d}(G)$ is the minimum number of color classes in a dominator coloring of $G$ (\cite{gera}). The dominator colorings were introduced by Gera \cite{gera} and studied by several authors. In particular, Kavitha and David \cite{kavitha} found dominator chromatic number of central graph of various graph families such as cycles, paths or wheel graphs. They also compared these parameters with  dominator chromatic number of the respective base graph families. Arumugam and Bagga \cite{Ar} obtained several result on graphs such that $\chi_d(G)= \chi(G)$ and $\chi_d(G)= \gamma(G)$. In particular, they proved that if $\mu(G)$ is the Mycielskian of G, then $\chi_d(G)+1 \leq \chi_d(\mu(G)) \leq \chi_d(G)+2$. Other results on the dominator coloring in graphs can be found in the works of Abdolghafurian et al. \cite{ab} (claw-free graphs), Bagan et al. \cite{g} (arbitrary graphs, $P_4$-sparse graphs, $P_5$-free graphs, claw-free graphs, graphs with bounded treewidth) Gera \cite{ge,ger} (various classes of graphs, in particular bipartite graphs), Merouane and Chellali \cite{h} (trees), Kazemi \cite{k} (total version of the problem) and Ramachandran et al. (e.g. windmill graphs). Askari et al. \cite{Ask} studied total global dominator coloring, where every vertex needs not only to dominate some color class other than its own, but also to be not adjacent to any vertex from another color class. Mojdeh et al. \cite{Moj} in turn considered strong dominator colorings, i.e., dominator colorings with an additional degree constraint, of central graphs and trees.

Instead of studying proper colorings, where all the neighbors of every vertex $v$ have color different than the color of $v$, one can focus on other kinds of coloring, e.g. majority coloring, where only half of the neighbors of every vertex need to have different color (or, which is equivalent in the case of finite graphs, at most half of the neighbors of any vertex $v$ can have the same color as $v$). It was proved long time ago by Lov\'asz that the majority coloring for any simple graph requires at most $2$ colors \cite{lovasz}. The problem was then considered for infinite graphs (see e.g. \cite{aharoni,shelah}) and recently for digraphs, where the out-neighbors were taken under consideration instead of all the neighbors, also in the list version, where every vertex can have a different list of colors. The best known results here say that four colors are enough for finite digraphs and five for infinite ones (see e.g. \cite{anholcer1, anholcer, seymour}).

Motivated by the above concepts, we decided to combine them. One could do this in several ways, by relaxing the condition of proper coloring or the condition of dominating the entire color class. Thus, one could require e.g. a majority coloring where every vertex dominates some color class. In this paper, however, we are interested in somehow opposite problem, where the coloring is supposed to be proper, but only half of some color class is supposed to be dominated by every vertex.

\begin{definition}\label{defMD}
\textit{Majority dominator coloring} of a graph $G$ is a proper coloring in which each vertex of $G$ dominates at least half of some color class. The majority dominator chromatic number $\chi_{md}(G)$ is the minimum number of color classes in a majority dominator coloring of $G$.
\end{definition}

Note that, according to the above definition, the domination condition is satisfied also if some vertex $v$ dominates its own color class. Of course, in such case this color class can consist of at most two vertices ($v$ and some $u\not\in N(v)$). 

The structure of the paper is as follows. In the next section we present some general results, in particular we characterize graphs $G$ for which $\chi_{md}(G)$ takes some specific values. In Section \ref{secFamilies} we derive the values of $\chi_{md}(G)$ for graphs $G$ from several families. We conclude the paper in the last section with some final remarks, conjectures and open problems. 

\section{General results}\label{secGeneral}

Since every majority dominator coloring is proper by definition, and on the other hand every dominator coloring is also a majority dominator coloring, we have that
$$
\chi(G)\leq \chi_{md}(G) \leq \chi_{d}(G) \leq n 
$$
for every graph $G$. Moreover, this bounds are sharp, as the example of the complete graph $K_n$ shows. The above implies in particular, that for every graph $G$ we have
$$
\chi_{md}\leq \chi(G) +\gamma(G),
$$
since the inequality $\chi_{d}\leq \chi(G) +\gamma(G)$ follows from \cite[Theorem 3.3]{ger}.

A well-known fact (see e.g. \cite[Lemma 3.1.33]{west}) is that a set of vertices  $S$ is an independent dominating set of $G$ if and only if it is a maximal independent set. This fact immediately implies that $\chi_{md}\leq \chi(G) +\alpha(G)$, but we can use it to prove a stronger inequality.

\begin{theorem}\label{m110}
Let $G$ be a connected graph. Then $\chi_{md}(G)\leq \chi(G) +\lceil \frac{\alpha}{2}\rceil-1$ and the bound is sharp.
\end{theorem}
\begin{proof}
Let $c:V(G)\rightarrow \{1,\dots,k\}$, where $k=\chi(G)$, be a proper $k$-coloring of $G$ and let $V_i=\{v\in V(G):c(v)=i\}$. Observe that there must be a proper coloring $c$ in which $V_k$ is maximal (not necessarily maximum) independent set (otherwise we can extend it by adding new vertices and recolor them with color $k$). Since $V_k$ is maximal independent, it is also dominating set in $G$. In particular, we have that $|V_k|\leq\alpha(G)$, $\chi(G[V(G)\setminus V_k])=\chi(G)-1$ (in fact, $G[V(G)\setminus V_k]$ has proper coloring with colors $1,\dots,k-1$ which is restriction of $c$) and every vertex from $V(G)\setminus V_k$ is adjacent to at least one vertex from $V_k$. Now, we leave two vertices from $V_k$ colored with $k$, and recolor the remaining $|V_k|-2$ of them arbitrarily with $\lceil(|V_k|-2)/2\rceil$ new colors so that each color appears at most twice. Since we are recoloring only vertices of one color class of $G$, the new coloring is a proper coloring with $k+\lceil(|V_k|-2)/2\rceil\leq\chi(G)+\lceil(\alpha(G))/2\rceil-1$ colors. Every vertex of $V_k$ has itself in its own closed neighborhood. Similarly, as we already observed, every vertex of $V(G)\setminus V_k$ has a neighbor in $V_k$. Since every color class in $V_k$ has at most two members, the obtained coloring is a majority dominator coloring. The sharpness of the bound can be observed e.g. for $G= C_n \circ K_1$, where $n\geq 4$ is even (i.e., $G$ consists of an even $n$-cycle, with $n$ extra pendant vertices, each of them adjacent to another cycle vertex). In these graphs $\chi(G)=2, \alpha(G)= n$, and $\chi_{md}(G)= \lceil n/2\rceil +1$ (however this property does not hold for odd $n$, see Proposition \ref{thmCirc} for details).
\end{proof}

Gera, Rasmussen and Horton \cite[Proposition 2.4]{gera} proved that for every connected $G$ of order $n\geq 3$, $\chi_d(G)\leq n+1-\alpha(G)$ holds, which implies in particular that $\chi_{md}(G)\leq n+1-\alpha(G)$. Also this inequality can be strengthened. Recall that $\overline{G}$ denotes the complement of $G$, that is the graph for which $V(\overline{G})=V(G)$ and $e\in E(\overline{G})\iff e\notin E(G)$.

\begin{theorem}\label{nMinusInd}
Let $G$ be a connected graph of order $n\geq 2$, where $I$ is an independent set in $G$ and $M$ is a matching in $\overline{G\setminus I}$. Then $\chi_{md}(G)\leq n-|M|-|I|+1$ and the bound is sharp.
\end{theorem}
\begin{proof}
Define the coloring $c:V(G)\rightarrow \{1,\dots,n-|M|-|I|+1\}$ as follows. Color all the vertices in $I$ with color $1$, the pairs of vertices in $M$ with distinct colors $2,\dots, |M|+1$ and the remaining $n-|I|-2|M|$ vertices of $V(G)\setminus (I\cup V(M))$ with distinct colors from the set $|M|+2,\dots, n-|I|-|M|+1$. This coloring is a majority dominator coloring, since it is proper, and every vertex in the graph dominates some vertex from $M$ (that is, half of the $2$-element color class) or a vertex in $V(G)\setminus (I\cup V(M))$ (that is, the entire $1$-element class). The sharpness of the bound can be observed e.g. for $G=K_n$, where $|I|=1$, $|M|=0$ and $\chi_{md}(G)=n$.
\end{proof}

This immediately implies the following.

\begin{corollary}\label{nMinusIndCol}
Let $G$ be a connected graph of order $n\geq 2$, where $I$ is a maximum independent set. Then
$\chi_{md}(G)\leq n +1-\alpha(G)-\nu(\overline{G\setminus I})$.
\end{corollary}

If $\Delta(G)=n-1$, then $\gamma(G)=1$ and by \cite[Lemma 3.4]{ger}, we get the following.

\begin{observation}\label{m112}
Let $G$ be a connected graph with $\Delta(G)= n-1$. Then
$$
\chi_{md}(G)=\chi(G).
$$
\end{observation}

We can show that the latter holds for a wider family of graphs.

\begin{proposition}\label{m113}
Let $G$ be a connected graph with $\Delta(G)\geq n-2$. Then $\chi_{md}(G) =\chi(G)$.
\end{proposition}
\begin{proof}
Since by Observation \ref{m112} we can assume that there is a vertex of degree $n-2$, say $v_1$, there is also exactly one vertex non-adjacent to $v_1$, say $v_2$. Consider any proper coloring $c$ of $G$ with $\chi(G)$ colors. The class of $c(v_1)$ has one or two members (in the latter case $c(v_1)=c(v_2)$). Define the coloring $c^\prime$ (possibly equal to $c$) as $c^\prime(v_2)=c(v_1)$ and $c^\prime(v)=c(v)$ for $v\neq v_2$. Obviously $c^\prime$ is proper, the class of $c^\prime(v_1)$ has exactly two members and every vertex $v\in V(G)\setminus\{v_1,v_2\}$ is adjacent to $v_1$,  thus $c^\prime$ is a majority dominator coloring.
\end{proof}

Let us conclude this section with a result for disconnected graphs.

\begin{proposition}
If $G$ be a disconnected graph with components $G_1, G_2, \ldots, G_k$ with $k\geq 2$ then\\
$$
\max_{j \in \{1, 2, \ldots, k\}} \chi_{md}(G_j)\leq \chi_{md}(G)\leq \sum_{j=1}^{k}(\chi_{md}(G_j)).
$$
The bounds are sharp.
\end{proposition}

\begin{proof}
For every $j$, $j \in \{1, 2, \ldots, k\}$, let $c_j:V(G_j)\rightarrow \{c_j^1,c_j^1+1,\dots,c_j^2=c_j^1+\chi_{md}(G_j)-1\}$ be a majority dominator coloring of component $G_j$, $j=1,\dots,k$. If the color sets of components are pairwise disjoint, we can assume that $c_j^1=c_{j-1}^2$ for $j=2,\dots,k$. Then obviously $\bigcup_{j=1}^{k} \{c_j\}$ is a majority dominator coloring of $G$ and $\chi_{md}(G)\leq \sum_{j=1}^{k}(\chi_{md}(G_j))$. To see that the bound is sharp, consider $G_i= P_4 \boxtimes P_2$ for $1\leq i \leq k$. Here we have $\chi_{md}(G_j)=4$ for every $j=1,\dots, k$ and none of the colors can be used in another $G_j$, since every vertex dominates exactly one color in every color class and there are exactly two vertices in every color class of $G_j$ (the coloring is unique up to obvious permutations of colors or vertices).

To prove the lower bound note the following. By definition, for any majority dominator coloring $c$ of $G$, the restriction of $c$ to $G_j$ is a majority dominator coloring of $G_j$, since the adjacency relations in any $G_j$ are exactly the same as in $G$. This implies that $\chi_{md}(G_j)\leq \chi_{md}(G)$ for every $j=1,\dots,k$.

To see that the bound is sharp, let $G_1= K_n$ for $n \geq 2$ and $G_j=K_2$ for every $j=2,\dots, n$. Here we have $\chi_{md}(G)=\chi_{md}(G_1)=n>2=\chi_{md}(G_j)$, $j=2,\dots,n$ (to see that $\chi_{md}(G)\leq n$, consider coloring $c:V(G)\rightarrow \{1,\dots,n\}$ in which for $j>1$, the vertices of $G_j$ are colored with $1$ and $j$.).
\end{proof}

Note that in order to prove the lower bound in the above theorem, we used the fact that the restriction of a majority coloring $c$ of $G$ to $G_j$ is a majority dominator coloring of $G_j$, because the adjacency relations are preserved. As one can see, it is not true for arbitrary subgraphs. To see it, note  that for example $\chi_{md}(C_{14})=5$, the colors of consecutive vertices being $1, 3, 2, 1, 3, 2, 1, 4, 2, 1, 4, 2, 1, 5$ (see Theorem \ref{thmPaths}, Lemma \ref{lem_propertiesCycle} and Corollary \ref{corCycles} for details). On the other hand, if we add the edge $v_7v_{14}$, then we obtain a majority dominator coloring with $4$ colors by recoloring $v_{14}$ with color $2$. This allows us to state the following.

\begin{observation}\label{subgraphsHer}
The majority dominator colorability is not a hereditary property, that is, it is not preserved by subgraphs in the general case. 
\end{observation}

In the reminder of this section we discuss the graphs $G$ such that $\chi_{md}(G)\leq 2$ or $\chi_{md}(G)\geq n-1$. Below we use $d_{S}(v)=\vert N(v)\cap S\vert$ to denote the number of neighbors of $v$ being members of a set $S$. Let us start with graphs $G$, for which $\chi_{md}(G)$ is close to its minimum possible value.

\begin{observation}\label{41}
Let $G$ be a graph. Then $\chi_{md}(G)= 1$ if and only if $G=\overline{K_n}$, $1\leq n\leq 2$.
\end{observation}
\begin{proof}
If $E(G)=\emptyset$, then obviously $\chi_{md}(G)=\left\lceil \frac {n}{2}\right\rceil$. If $E(G)\neq\emptyset$, then $\chi_{md}(G)\geq \chi(G)\geq2$.
\end{proof}

\begin{proposition}
A graph $G\notin\{{K_1}, \overline{K_2}\}$ satisfies $\chi_{md}(G)= 2$ if and only if it is bipartite, where $X$ and $Y$ are partition sets such that at least one of the following holds:
\begin{enumerate}[label=(\roman*)]
\item
$1\leq |X| \leq |Y| \leq 2$,
\item
$1\leq |X| \leq 2<|Y|$ and for every $y\in Y$, $d(y)\geq 1$, 
\item
$3\leq |X|\leq |Y|$, and for every $x\in X$ and $y\in Y$, $d(x)\geq |Y|/2$ and $d(y)\geq |X|/2$, respectively.
\end{enumerate}
\end{proposition}
\begin{proof}
By Observation \ref{41}, for every graph $G$ from the described families we have $\chi_{md}(G)\geq 2$. On the other hand, it is a routine to check that in every case the coloring assigning color $1$ to the vertices of $X$ and color $2$ to the vertices of $Y$ is a majority dominator coloring and thus $\chi_{md}(G)\leq 2$.

Now assume that $\chi_{md}(G)=2$. Obviously $G\notin\{{K_1}, \overline{K_2}\}$ by Observation \ref{41} and $G$ is bipartite since $\chi(G)\leq 2$. Let $X$ and $Y$ be the partition sets. Each vertex must dominate half of the vertices in at least one partition set: its own (but then the order of this set must be at most $2$) or the other one (with at least half of the vertices of this set being its neighbors in this case). Thus, three cases are possible: $|X|\leq |Y|\leq 2$ (covered by (i)), $|X|\leq 2$ and $|Y|\geq 3$, while every vertex in $Y$ dominates at least half of the vertices in $X$ (covered by (ii)) and both $|X|,|Y|\geq 3$, while every vertex in $X$ (and $Y$) dominates at least half of the vertices in $Y$ ($X$, respectively).
\end{proof}

The following is immediate.

\begin{corollary}
For every $1\leq m\leq n$, we have $\chi_{md}(K_{m,n})=2$.
\end{corollary}

Now we switch to the graphs $G$, for which $\chi_{md}(G)$ takes maximum possible values.

\begin{proposition}\label{chiN}
A graph $G$ satisfies $\chi_{md}(G)= n$ if and only if $G=K_n$.\\
\end{proposition}

\begin{proof}
The only proper coloring of $K_n$ is also its majority dominator coloring, thus $\chi_{md}(K_n)=\chi(K_n)=n$. On the other hand, assume that for some graph $G$, every vertex is colored with a different color. If some two vertices $x$ and $y$ are not adjacent, one can recolor $y$ with the color of $x$ and obtain this way a majority dominator coloring (where each vertex dominates at least half of its own color class) with $n-1$ colors, thus $\chi_{md}(G)\leq n-1$. This means that if $\chi_{md}(G)=n$, then there are no non-adjacent vertices and so $G=K_n$.
\end{proof}

Let $E(v)$ denote the set of edges incident with given vertex $v\in G$.

\begin{proposition}\label{chiN1}
Let $G$ be a connected graph of order $n\geq 3$. Then $\chi_{md}(G)= n-1$ if and only if $G=K_n-E^\prime(v)$ for some set $E^\prime(v)\subset E(v), v\in V(G)$ such that $1\leq |E^\prime(v)|\leq n-2$.
\end{proposition}

\begin{proof}
Let $G=K_n-E^\prime(v)$. Observe that $G$ is exactly $K_{n-1}$ with one extra $v$, where $1\leq d(v)\leq n-2$. It means that $\Delta(G)=n-1$ and $\chi(G)=n-1$, so by Proposition \ref{m113}, $\chi_{md}(G)=n-1$.

Now, Assume that there is a connected graph $G$ such that $\chi_{md}(G)= n-1$. We claim that $\Delta(G)=n-1$. Indeed, if $d(u)\leq n-2$ for every $u\in V(G)$, then there are at least two non-adjacent vertices $u_1, u_2$ of degree at most $n-2$ in $G$. If $G\setminus\{u_1,u_2\}$ is not a complete graph, then by Proposition \ref{chiN} there exists a majority dominator $(n-3)$-coloring $c$ of it and we can extend it to a majority dominator $(n-2)$-coloring $c^\prime$ of $G$ by putting $c^\prime (u_1)=c^\prime (u_2)=n-2$, a contradiction. On the other hand, if $G\setminus\{u_1,u_2\}=K_{n-2}$, then $N(u_1)\cap N(u_2)=\emptyset$, because $d(v)\leq n-2$ implies $d_{\{u_1,u_2\}}(v)\leq 1$ for every $v\in V(G\setminus\{u_1,u_2\})$. For any (majority dominator) $(n-2)$-coloring $c$ of $G\setminus\{u_1,u_2\}$, and for any $v_1\in N(u_1)$, $v_2\in N(u_2)$ we can extend it to a majority dominator $(n-2)$-coloring $c^\prime$ of $G$ by setting $c^\prime(u_1)=c(v_1)$ and $c^\prime (u_2)=c(v_2)$, a contradiction. By Observation \ref{m112}, $\Delta(G) =n-1$ implies $\chi(G)=\chi_{md}(G)=n-1$.

If there are two disjoint edges $u_1v_1,u_2v_2\in E(\overline{G})$, then $\chi(G)\leq n-2$, since one can set $c(u_1)=c(v_1)\neq c(u_2)=c(v_2)$ and color the remaining $n-4$ vertices with other $n-4$ colors. Also, there is no triangle in $\overline{G}$ since in such a case all three of its vertices could get the same color and we would obtain again $\chi(G)\leq n-2$. Thus the edges in $\overline{G}$ must form a star having at least one and at most $n-2$ edges. This ends the proof.

\end{proof}


\section{Majority dominator coloring of chosen graphs}\label{secFamilies}

In this section we are going to present the exact values of $\chi_{md}(G)$ for chosen graphs $G$. From the previous section we already know that, in particular
$\chi_{md}(K_n)=n$ and $\chi_{md}(K_{m,n})=2$. We also referred to the following fact.

\begin{proposition}\label{thmCirc}
For $n\geq 3$, $\chi_{md}(C_n\circ K_1)= \left\lceil\frac{n}{2}\right\rceil+1$.
\end{proposition}
\begin{proof}
Let $n\geq 3$ and let $G=C_n\circ K_1$. From Theorem \ref{m110} it comes that for even $n$, $\chi_{md}(G)\leq \lceil n/2\rceil+1$, since $\chi(G)=2$ and $\alpha(G)=n$. When $n$ is odd, we have $\chi(G)=3$ and from Theorem \ref{m110} it follows that $\chi_{md}(G)\leq \lceil n/2\rceil+2$, but one can reduce this bound by $1$. Denote the consecutive vertices of $C_n$ by $v_1, \dots, v_n$ (so that $v_iv_{i+1}\in E(G)$ for $i=1,\dots,n-1$ and $v_1v_n\in E(G)$) and their neighbors by $u_1,\dots, u_n$ (so that $v_iu_i\in E(G)$ for $i=1,\dots, n$). Let us define the coloring $c$ as: $c(v_{2i+1})=c(u_{2i})=i$ for $i=1,\dots, (n-1)/2$, $c(v_1)=c(u_n)=(n+1)/2$ and $c(v_{2i})=c(u_{2i-1})=(n+3)/2=\lceil n/2\rceil+1$ for $i=1,\dots,n-1$. Clearly, $c$ is an $(\lceil n/2\rceil+1)$-majority dominator coloring of $G$.

Now we are going to show that at least $\lceil n/2\rceil+1$ colors are necessary for every $n$. Consider the pendant vertices. Each of them can dominate at most $1$ vertex -- either itself or its only neighbor. In both cases it comes out that at most $2$ vertices can be colored with the dominated color. It means that at least $\lceil n/2\rceil$ colors must be dominated by the pendant vertices and that at most $2\lceil n/2\rceil\leq n+1$ vertices can be colored with these colors. The remaining $2n-2\lceil n/2\rceil\geq n-1\geq 2$ vertices need at least one extra color.
\end{proof}

In the case of wheel $W_n$, every proper coloring is also a majority dominator coloring (in particular, each vertex dominates the one-element class of the central vertex). This implies that for every $n\geq 3$, $\chi_{md}(W_n)= \chi(W_n)$.

\begin{observation}\label{lem21}
For $n\geq 3$,
$$
\chi_{md}(W_n)=
\begin{cases}
4, &\text{if } n\equiv 1 \imod 2,\\
3, &\text{if } n\equiv 0 \imod 2.
\end{cases}
$$
\end{observation}

Let $S_{a,b}$ be the double star with central vertices $u$ and $v$ with $d(u)=a\geq 2$ and $d(v)=b\geq 2$. Let $X=\{x_1, x_2, \dots, x_{a-1}\}$, $Y=\{y_1, y_2, \dots, y_{b-1}\}$, $N(v)= \{u\} \cup Y$ and $N(u)=\{v\}\cup X$.
\begin{proposition}
For $a,b\geq 3$, $\chi_{md}(S_{a,b})=3$. Otherwise, $\chi_{md}(S_{a,b})=2$.
\end{proposition}

\begin{proof}
If $a=2$ or $b=2$, then the only proper $2$-coloring defined by $V_1= X\cup \{v\}$ and $V_2=Y\cup \{u\}$ is a majority dominator coloring (at least one of the color classes consists of two members and at least one of them is dominated by every vertex in the graph). On the other hand, if $a,b\geq 3$, then $u\in V_1$ and $v\in V_2$. Assume that these two color classes are enough. Then $V_1=\{u\}\cup Y$ and $V_2=\{v\}\cup X$. We have $|V_1|,|V_2|\geq 3$. But every vertex in $X\cup Y$ can dominate only one vertex and thus a color class of at most two members, a contradiction. It follows that at least three colors are necessary. On the other hand, the coloring defined by $V_1=\{u\}$, $V_2=\{v\}$, $V_3=X\cup Y$ is a majority dominator $3$-coloring of $S_{a,b}$.
\end{proof}

The star $K_{1, n-1}$ has one vertex $v$ of degree $n-1$ and $n-1$ vertices of degree one. It can be generalized to the multistar graph $K_{n}(a_{1}, a_{2}, . . ., a_{n})$, which is formed by attaching $a_{i} \geq 1$, ($1 \leq i \leq n$) pendant vertices to each vertex $x_{i}$ of a complete graph $K_{n}$, with $V(K_{n})=\{x_{1}, x_{2}, . . ., x_{n}\}$.

\begin{proposition}
For the multistar graph $G= K_{n}(a_{1}, a_{2}, . . ., a_{n})$\\
{
 $$
\chi_{md}(G)=
\begin{cases}
n, &\text{if } a_i<n \text{ for some } 1\leq i\leq n,\\
n+1, &\text{otherwise}.
\end{cases}
$$
}
\end{proposition}

\begin{proof}
It is clear that $\chi(G)=n$, so by $\chi(G) \leq \chi_{md}(G)$ it follows that $\chi_{md}(G) \geq n$. Let $V(K_{n})=\{x_{1}, x_{2}, . . ., x_{n}\}$. We will use $A_i$ to denote the set of pendant vertices adjacent to $x_i$. Note that $a_i=|A_i|$.

Assume first that there is $a_i$ such that $a_i<n$ (without loss of generality, assume that $i=n$). Color $K_n$ properly with $n$ colors so that $c(x_j)=c_j, j=1,\dots, n$, then all the vertices in $A_j, j<n$ with $c_n$, and finally the vertices of $A_n$ with different colors $c_1,\dots,c_{a_n}$. Every vertex in $G$ dominates at least one vertex colored with some $c_j, j<n$, that is at least half of the respective color class, so $c$ is a majority dominator coloring.

Assume now that for each  $A_i$ we have $a_i \geq n$. Suppose that  $\chi_{md}(G) = n$. Obviously, we have $c(x_i)=c_i$ for $i=1,\dots n$. If any of these colors, say $c_j$, appears at least $3$ times in $G$, then all the pendant neighbors of $x_j$ must belong to $2$-color classes, since they cannot dominate the class of $x_j$. However, each of those colors appears already in $K_n$, so there would be at least $n$ different colors in $A_j$ different than $c_j$, a contradiction. Thus at least $n+1$ colors are necessary. On the other hand the coloring defined as $c(x_i)=c_i$ for $i=1,\dots n$ and $c(v)=c_{n+1}$ for $v\in \bigcup_{i=1}^n{A_i}$ is a majority dominator $(n+1)$-coloring.
\end{proof}

In the reminder of this section we are going to focus on the majority dominator coloring numbers of paths and cycles. Let as start with the following technical result.

\begin{lemma}\label{lem_propertiesPath}
Let $P=P_n$ be a path with vertex set $V(P)=\{v_1, v_2, \dots, v_n\}$, let $k=\chi_{md}(P)$ and let $c:V(P)\rightarrow\{1,\dots,k\}$ be its $k$-majority dominator coloring. Then the following statements are true.
\begin{enumerate}
    \item 
    There are at most two colors in $c$ that appear on at least five vertices.
    \item
    If $n\geq 11$ then there are at least two colors in $c$ that appear at most twice.
    \item
    If $n\neq 2$ then there is at most one color in $c$ that appears only once.
\end{enumerate}
\end{lemma}

\begin{proof}
In all three parts of the proof we are going to use the fact that $k$ is the minimum number of colors that can be used in a majority dominator coloring of $P$.\\

Part 1.\\
If $n<15$ then the conclusion is trivial. If $n\geq 15$, suppose that there are three colors $c_1$, $c_2$ and $c_3$ such that they occur at least five times in $c$. Since every vertex has degree at most $2$ and so it can dominate a color class of order at most $4$, it is impossible that there are three consecutive vertices $v_i, v_{i+1}, v_{i+2}$ such that $c(\{v_i, v_{i+1}, v_{i+2}\})\subseteq \{c_1,c_2,c_3\}$, because in such a case $v_{i+1}$ could not dominate any color class. For that reason the vertices colored with $c_1$, $c_2$ and $c_3$ can have at most one neighbor colored with one of those colors. It means that one can recolor every vertex colored with $c_3$ using $c_1$ (by default) or $c_2$ if one of its neighbors is colored with $c_1$. The new coloring $c^\prime$ is obviously proper. It is also a majority dominator coloring, since every vertex in $V(P)$ had to dominate a color class from $\{1,\dots,k\}\setminus\{c_1,c_2,c_3\}$ before the recoloring. Moreover it uses $k-1$ colors, which contradicts with the minimality of $c$ and concludes the proof.\\

Part 2.\\
Since $v_1$ and $v_n$ can dominate only a color class of order at most two, there must be at least one such color class in $c$. Suppose that there is only one such class, colored with color $c_1$. Since $v_1$ and $v_n$ both need to dominate this class, we must have in particular that $c_1\in c(\{v_1,v_2\})$ and $c_1\in c(\{v_{n-1},v_n\})$. No matter which of $v_1$ and $v_2$ is colored with $c_1$, the color of the other one, say $c_2$, must appear at least three times. Thus there must be at least one vertex $v_i$ colored with $c_2$, $3\leq i\leq n-2$.

If $n\geq 12$ then two cases are possible. In the first case there exists $v_i$ colored with $c_2$ having a neighbor $v_j$, $j\in\{i-1,i+1\}$ such that $c(v_j)=c_3$, $c_3\not\in \{c_1,c_2\}$, and $v_j$ has exactly one neighbor colored with $c_2$ and no neighbor colored with $c_1$, so also exactly one neighbor colored with yet another color $c_4$. This means that $v_j$ needs to dominate another color class of order at most two that must be $c_3$ or $c_4$, a contradiction. In the second case every neighbor of a vertex colored with $c_2$ has two neighbors of this color, but this means that there are at least $4$ vertices colored with $c_2$, $3\leq i\leq n-2$ and so at least five such vertices in total, thus any of their neighbors has to dominate its own class which must be of order at most two, a contradiction.

If $n= 11$ then there could be only two ways to avoid the multiple occurrence of a color class of order at most $2$. Either there are only three colors in use or four colors are used and all the color classes but one have cardinalities $3$. In the first case, however, ther must be a vertex such that its color (say $c_2$) appears at least three times and the color of its both neighbors (say $c_3$) appears at least five times, so it cannot dominate any color class. In the other case there is a vertex colored $c_2$ having two neighbors colored $c_3$ and $c_4$, where each color class has exactly three members, so it does not dominate any color class.\\

Part 3.\\
If $n\leq 4$ then the statement follows from the fact that every optimal proper coloring is majority dominator coloring. Thus we may assume that $n\geq 5$. Suppose that there are two colors $c_1$ and $c_2$ such that $|c^{-1}(c_1)|=|c^{-1}(c_2)|=1$. Let $v_i$ and $v_j$, $1\leq i<j\leq n$ be the two vertices such that $c(v_i)=c_1$ and $c(v_j)=c_2$. If they are not adjacent, then we can set $c(v_j)=c_1$ and obtain this way a new majority dominator coloring with $k-1$ colors, a contradiction. So assume that $j=i+1$. Since $n\geq 5$, without loss of generality we can assume that $i\geq 3$. We define a new coloring $c^\prime$ by $c^\prime(v_{i-1})=c_1$, $c^\prime(v_{i})=c(v_{i-1})$, $c^\prime(v_{i+1})=c_1$ and $c^\prime(v)=c(v)$ for $v\in V(P)\setminus \{v_{i-1},v_{i},v_{i+1}\}$. Note that the new coloring is a majority dominator coloring, since the vertices $v_{s}$, $i-2\leq s\leq \min\{i+2,n\}$ dominate now the color class $c_1$ while the other vertices (if there are any) dominate the same classes as before the recoloring. Moreover the new coloring uses $k-1$ colors, which contradicts with the minimality of $c$. Thus the number of one-member color classes cannot be greater than one.
\end{proof}

\begin{lemma}\label{lemmaPminus6}
Let $P=P_{n}$ be a path of order $n\geq 11$ with vertices $v_1, v_2, \ldots, v_{n-1}, v_{n}$, let $k=\chi_{md}(P)$ and let $c:V(P)\rightarrow\{1,\dots,k\}$ be its $k$-majority dominator coloring. Let $P^\prime=P_{n+6}$ be an extension of $P$ with vertices $v_1, v_2, \ldots, v_{n+5}, v_{n+6}$. Then $\chi_{md}(P^\prime) \geq k+1$.
\end{lemma}
\begin{proof}
Suppose that there exists a $k$-majority dominator coloring $c^\prime$ of $P^\prime$. We will show that this would contradict with the minimality of $c$. We will consider two cases: when $v_{n+6}$ dominates the color class of $c(v_{n+6})$ and when it does not.\\

Case 1. $v_{n+6}$ dominates its own color class.\\
Since all the colors used in $c^\prime$ must occur in its restriction to $P$ (otherwise it would use less than $k$ colors), the color $c_1=c^\prime(v_{n+6})$ appears exactly once in $P_n$. By Lemma \ref{lem_propertiesPath} it is the only one such color, thus every other color used in $V(P^\prime)\setminus V(P)$ must be used at least two times in $P$, so at least three times in $P^\prime$. This means in particular that each vertex $v_i$, $n+2\leq n+4$ must dominate a color class with at least three members and thus each of them must have two neighbors with the same color. In particular $c(v_{n+3})=c(v_{n+5})$ (because of $v_{n+4}$) and $c(v_{n+1})=c(v_{n+3})$ (because of $v_{n+2}$). This means that the color $c_2=c(v_{n+1})$ is used three times in $V(P^\prime)\setminus V(P)$ and at least two times in $V(P)$, so in fact it cannot be dominated, a contradiction.\\

Case 2. $v_{n+6}$ does not dominate its own color class.\\
In such situation $v_{n+6}$ has to dominate the color class of $v_{n+5}$. Using similar argument as in case $1$ we conclude that the vertices in $V(P^\prime)\setminus V(P)$ have to be colored in the following way: $c(v_{n+6})=c_1$, $c(v_{n+5})=c_2$, $c(v_{n+2})=c(v_{n+4})=c_3$ and $c(v_{n+1})=c(v_{n+3})=c_4$, where $c_1$, $c_2$, $c_3$ and $c_4$ are pairwise distinct. In particular $c_1\neq c_3$, because otherwise $c_1$ would appear three times in $V(P^\prime)\setminus V(P)$ and according to Lemma \ref{lem_propertiesPath} at least two times in $V(P)$, so $v_{n+3}$ could not dominate any color class. Also, in the restriction of $c^\prime$ to $P$ $c_1$ appears at least two times, $c_2$ exactly one time and $c_3$ and $c_4$ exactly two times. We are going to show that we can recolor the vertices colored with $c_1$, $c_3$ or $c_4$ in $P$ and obtain this way a $(k-1)$-majority dominator coloring of $P$, which contradicts the minimality of $c$.

We will consider two subcases: when some vertex dominates the color class of $c_1$ in $c^\prime$ and when there is no such vertex.\\

Subcase 2.1. There is a vertex $v\in V(P)$ that dominates the color class $c_1$ in $c^\prime$.\\
In such a case $c_1$ appears at most three times in $V(P)$ and moreover there exists a sequence of consecutive vertices $u,v,w$ in $V(P)$ such that $c(u)=c(w)=c_1$. Depending on the neighborhood of $u$ and $w$, we can recolor each of these vertices either to $c_3$ or to $c_4$ and obtain this way a new majority dominator coloring $c^{\prime\prime}$ of $P$. Indeed, if any vertex dominated the color class $c_3$ or $c_4$ in $c^\prime$, then it will also dominate it in $c^{\prime\prime}$, since after the recoloring we have $|(c^{\prime\prime})^{-1}(c_3)|\leq 4$ and $|(c^{\prime\prime})^{-1}(c_4)|\leq 4$. Now, either $c^{\prime\prime}$ is a $(k-1)$-coloring (if $c_1$ appears exactly two times in the restriction of $c^\prime$ to $V(P)$) or there is still one vertex colored $c_1$, which contradicts Lemma \ref{lem_propertiesPath}. Note that the recoloring step is always possible. The only problem can occur if both $u_1$ and $u_2$ have neighbors in both colors $c_3$ and $c_4$, but it can happen only if there exists a sequence of consecutive vertices $v_1,u_1,v,u_2,v_2$ colored with $c_4,c_1,c_3,c_1,c_4$ or $c_3,c_1,c_4,c_1,c_3$, respectively. However, such a configuration is impossible since in such a case $u_1$ would not dominate any color class in $c^\prime$.\\

Subcase 2.2. There is no vertex dominating the color class $c_1$ in $c^\prime$.\\
If there is no vertex $v\in V(P)$ that dominates $c_3$ in $c^\prime$, then we can obtain a new coloring $c^{\prime\prime}$ of $V(P)$ by recoloring both vertices in $(c^\prime)^{-1}(c_3)$ with $c_1$. This will not change any domination relation, so $c^{\prime\prime}$ is a $(k-1)$-majority dominator coloring of $P$, a contradiction. The same reasoning applies if we substitute $c_3$ with $c_4$. Thus the only case to consider is the situation when there is a sequence of consecutive vertices in $V(P)$, $u_1,v_1,w_1,S,u_2,v_2,w_2$, where $S$ denotes some subsequence of vertices (possibly empty), $c^\prime(u_1)=c^\prime(w_1)=c_3$ and $c^\prime(u_2)=c^\prime(w_2)=c_4$ (or $c^\prime(u_1)=c^\prime(w_1)=c_4$ and $c^\prime(u_2)=c^\prime(w_2)=c_3$, but we can skip this case by symmetry). If $S$ is not empty, then we can define a new $(k-1)$-majority dominator coloring $c^{\prime\prime}$ of $V(P)$ by recoloring $c^{\prime\prime}(w_1)=c^{\prime\prime}(w_2)=c_3$, a contradiction. If $S$ is empty, then we have a sequence $u_1,v_1,w_1,u_2,v_2,w_2$ such that $c^\prime(u_1)=c^\prime(w_1)=c_3$ and $c^\prime(u_2)=c^\prime(w_2)=c_4$, where $w_1$ dominates $c^\prime(v_1)$, so $|(c^\prime)^{-1}(c^\prime(v_1))|\leq 2$ and $c^\prime(v_1)$ is also dominated by $u_1$. By symmetry, $|(c^\prime)^{-1}(c^\prime(v_2))|\leq 2$ and $c^\prime(v_2)$ is dominated by $u_2$ and $w_2$. At least one of $c^\prime(v_1)$, $c^\prime(v_2)$ is not $c_2$. Without loss of generality assume that $c^\prime(v_1)\neq c_2$. Now, define the new coloring $c^{\prime\prime}$ of $V(P)$ by setting $c^{\prime\prime}(u_1)=c_1$ and $c^{\prime\prime}(w_1)=c_2$ (or $c^{\prime\prime}(u_1)=c_2$ and $c^{\prime\prime}(w_1)=c_1$ if $u_1$ has a neighbor colored with $c_1$). In both cases, after the recoloring, $u_1$, $v_1$ and $w_1$ dominate the color class of $c^{\prime\prime}(v_1)$ and other domination relations remain intact. It follows that $C^{\prime\prime}$ is a $(k-1)$-majority dominator coloring of $P$, a contradiction.
\end{proof}

\begin{theorem}\label{thmPaths}
Let $P=P_n$ be a path of order $n$ with vertices $v_1, v_2, \ldots, v_{n-1}, v_{n}$. Then
$$
\chi_{md}(P)=
\begin{cases}
1,& \text{if }n=1,\\
2,& \text{if }2\leq n\leq 5,\\
3,& \text{if }6\leq n\leq 10,\\
4,& \text{if }11\leq n\leq 13,\\
\left\lceil\frac{n}{6}\right\rceil+2,& \text{if }n\geq 14.
\end{cases}
$$
\end{theorem}

\begin{proof}
First we are going to show that the given values are minimum possible. For $n\leq 5$ it is obvious, since every coloring must be proper. On the other hand it is impossible to use only two colors for $n\geq 6$, since both colors would appear at least three times and the ends of the path have only one neighbor, so they cannot dominate any color. Assume that three colors are enough for $n=11$. By Lemma \ref{lem_propertiesPath} at least two of them occur at most two times, so the third one needs to be used at least $7$ times. However this would force at least one pair of neighboring vertices to be colored with the same (third) color, so one needs at least four colors in this case (as well as for every $n\geq 11$).

Assume that in a coloring $c$ of $P_{14}$ one can use four colors. By Lemma \ref{lem_propertiesPath} there are two or three colors that occur at most $2$ times. In the latter case, however, there would be at least eight vertices colored with the same (fourth) color, which would force at least one pair of them to be adjacent. Thus there are exactly two colors (say $c_1$ and $c_2$) present on at most two vertices and at least ten vertices colored with two other colors (say $c_3$ and $c_4$). Since in the sequence of colors there can be at most five subsequences consisting of $c_3$ and $c_4$ only (separated by at most four vertices colored $c_1$ or $c_2$) either all of them are pairs, or there are some consisting of at least three vertices. The former is impossible, since then each of $c_3$ and $c_4$ occurs exactly five times and the ends of the path (colored with $c_3$ or $c_4$ and neighboring with vertices colored with $c_4$ or $c_3$, respectively) cannot dominate any color class. Thus there must be at least one subsequence consisting of at least three vertices. Without loss of generality assume that it is $c_3, c_4, c_3$. Because of the middle vertex, $|c^{-1}(c_3)|\leq 4$, so $|c^{-1}(c_4)|\geq 6$. But it means that in at least one of the five subsequences $c_4$ occurs at least twice, thus there is a subsequence of colors $c_4, c_3, c_4$ and we deduce that $|c^{-1}(c_4)|\leq 4$, a contradiction. This means that at least five colors are necessary for every $n\geq 14$.

Finally, assume that in a coloring $c$ of $P_{19}$ one can use five colors. According to Lemma \ref{lem_propertiesPath} the number of colors present on at most two vertices is between two and four. If it is four, then there are at least eleven vertices colored with the fifth color, which forces at least two of them to be adjacent, a contradiction. If there are three such colors, then at least thirteen vertices must be colored with two other colors (say $c_1$ and $c_2$) and form at most seven subsequences. If they are distributed in six pairs and one singleton, then there is at least one end of the path that does not dominate any color. Thus there must be at least one triple, say $c_1, c_2, c_1$. This implies that there are at most four vertices labeled with $c_1$ and consequently at least nine vertices colored with $c_2$. Thus at least one of (at most) seven subsequences contains at least two of them, so there is a subsequence $c_2, c_1, c_2$, which implies that $c_2$ occurs at most four times, a contradiction. Finally, assume that there are exactly two colors occurring at most two times. Then there are at least fifteen vertices colored with three colors, say $c_1$, $c_2$ and $c_3$. By Lemma \ref{lem_propertiesPath} there are at most two colors appearing on at least five vertices, so without loss of generality we can assume that $c_1$ occurs at most four times and $c_2$ at least six times. Observe that a sequence being a permutation of $c_1,c_2,c_3$ is impossible since the middle vertex could not dominate any color class. This implies that there is a subsequence $c_2,x,c_2$, where $x\in\{c_1,c_3\}$. It implies that $c_2$ can occur at most $4$ times, a contradiction. We deduce that at least six colors are necessary in any majority dominator coloring of $P_n$, where $n\geq 19$.

As we can see, in particular we obtained that $\chi_{md}(P_n)\geq \left\lceil\frac{n}{6}\right\rceil+2$ for $14\leq n\leq 19$. Using Lemma \ref{lemmaPminus6} it follows by induction that $\chi_{md}(P_n)\geq \left\lceil\frac{n}{6}\right\rceil+2$ for every $n\geq 14$.

To conclude the proof, we need to show that there exist majority dominator colorings of $P_n$ using the number of colors given in the statement of the Theorem.

As it can be easily verified, for $n\leq 13$ we can use the sequences of colors $S_n$ defined as follows: $S_1=\{1\}$, $S_2=\{1,2\}$, $S_3=\{1,2,1\}$, $S_4=\{1,2,1,2\}$, $S_5=\{1,2,1,2,1\}$, $S_6=\{1,2,1,2,1,3\}$, $S_7=\{3,1,2,1,2,1,3\}$, $S_8=\{3,1,2,1,2,1,2,3\}$, $S_9=\{3,1,2,1,2,1,2,3,1\}$, \\$S_{10}=\{3,1,2,1,2,1,2,1,3,2\}$, $S_{11}=\{3,4,1,2,1,2,1,2,1,2,3\}$, \\$S_{12}=\{3,4,1,2,1,2,1,2,1,2,3,4\}$, $S_{13}=\{1,2,1,2,1,3,1,3,1,4,1,4,1\}$.

For $n\geq 14$ the desired coloring can be defined with the following formula:
$$
c(v_i)=
\begin{cases}
1,&\text{if } i\equiv 0\imod 3,\\
2,&\text{if } i\equiv 1\imod 3 \wedge i\neq n,\\
\left\lceil\frac{n}{6}\right\rceil+2,&\text{if } i\equiv 1\imod 3 \wedge i=n,\\
\left\lceil\frac{i}{6}\right\rceil+2,&\text{if } i\equiv 2\imod 3.
\end{cases}
$$
Note that according the above formula for every $1\leq j\leq n$ we have $c(v_{j-1})\geq 3$, $c(v_{j})\geq 3$ or $c(v_{j+1})\geq 3$ and each color $x\geq 3$ appears at most two times, so $c$ is indeed a majority dominator coloring of $P$.
\end{proof}

In order to present the result for cycles, we will need the following technical lemma.

\begin{lemma}\label{lem_propertiesCycle}
Let $C=C_n$, $n\geq 9$ be a cycle with vertex set $V(C)=\{v_1, v_2, \dots, v_n\}$, let $k=\chi_{md}(C)$ and let $c:V(C)\rightarrow\{1,\dots,k\}$ be its $k$-majority dominator coloring. Let $P=P_n$ be a path of order $n$. Then $k\geq \chi_{md}(P_n)$.
\end{lemma}

\begin{proof}
Obviously there must be a color that appears once or twice in $c$. Otherwise every vertex would dominate a color present on at least three vertices being the color of its both neighbors but this could occur only if the cycle were colored alternately with two colors, which is possible only for $n\in\{4,6,8\}$.

If some color occurs exactly once, say on vertex $v$, then we are done, since we can remove any edge incident with $v$, say $uv$, and recolor $u$ with $c(v)$, which defines a majority dominator coloring of $P_n$.

From now on we will assume that there are no colors appearing only once. First assume that there exist two vertices $v_1$ and $v_2$ such that $|c^{-1}(c(v_1)|\leq 2$, $|c^{-1}(c(v_2)|\leq 2$ and distance between $v_1$ and $v_2$ is at most $3$. In such a case it is enough to remove an edge from a shortest $(v_1,v_2)$-path: either the edge incident to both $v_1$ and $v_2$ (if the distance is $1$) or an edge not incident to at least one of $v_1$ and $v_2$. This way each of $v_1$ and $v_2$ will become either an end vertex of the obtained path or the only neighbor of an end vertex and $c$ will still be a majority dominator coloring of the obtained path.

It remains to focus on the case where every color that appears on at most two vertices, appears on exactly two vertices, which means that there are exactly $2p$ vertices colored with colors appearing twice, where $p$ is the number of such colors, and the distance between each pair of such vertices is at least four, or in other words there is a sequence of at least three consecutive vertices colored with some colors used at least three times between every such pair (we will call such sequences \textit{segments}). Note that the number of segments is also $p$. It is impossible that three colors $c_1$, $c_2$, $c_3$ are used on three consecutive vertices, because in such case the middle vertex would not dominate any color class. This implies in turn that in every segment exactly two colors are used. If it has the form $c_1,c_2,c_1$, then the color $c_1$ can appear at most $4$ times in $c$. If there are at least four colors in a segment, then each od $c_1$ and $c_2$ can appear at most four times (so in particular the length of the segment cannot be more than eight).

If every segment has length $3$, then one can use one common color on the middle vertex of every segment and one color on the remaining vertices of every pair of segments. This allows to use $p+1$ colors in the segments and $k=2p+1$ colors in total to color $n=8p$ vertices. Since exactly two colors must appear in every segment, this is the minimum number of colors that can be used. Note that $k=\frac{n}{4}+1$ in this case.

Assume that there is a segment of length at least $4$. Then the number of segments of length three equals to $s_3\in\{2q,2q+1\}$ for some integer $q<p$. Now, let us consider a few cases.
\begin{enumerate}
    \item The lengths of the remaining $2p-s_3$ segments are all $4$ or there are possibly $s_5\leq s_3$ segments of length $5$. We color optimally $2\left\lfloor\frac{s_3-s_5}{2}\right\rfloor=2q^\prime$ segments of length $3$ with $q^\prime+1$ colors, like above. To color the remaining segments we need $2(p-q^\prime)$ new colors (two for every pair, where both segments have length $4$ or one of them has length $3$ and another $5$ or $4$, where the last case can occur at most once, when $s_3-s_5$ is odd), so in total $k=p+(q^\prime+1)+2(p-q^\prime)=3p-q^\prime+1$ colors are necessary to color $n\leq 2p+3\times 2q^\prime+4\times 2(p-q^\prime)=10p-2q^\prime$ vertices. Note that since $q^\prime\leq q\leq p-1$, we have $k\geq \frac{5}{2}p-\frac{1}{2}q^\prime+\frac{3}{2}\geq \frac{n}{4}+1$ in this case.
    \item There are some segments of at least one of the following kinds: segments of length $6$ or more (let us denote their number by $s_6$), possibly some segments of length $5$ that could not be paired in the process described above when $s_5>s_3$ (in such case let $s_6^\prime=s_6+s_5-s_3$) or finally possibly one unpaired segment of length $4$ (in such case increase $s_6^\prime$ by $1$). Note that each of these $s_6^\prime$ segments not existing in the previous cases needs two new colors, with the only exception when there is one unpaired segment of length $4$ and at least one unpaired segment of length $5$, when they can be paired and colored with three new colors. It follows from the fact that if a segment has length at least $5$, then the two colors used in it cannot be reused to color any segment of length at least $4$ with the only exception mentioned above, when one of them can be used again. Since the length of any segment cannot be more than $8$, one can easily see that the total number of colors $k\geq p+(q^\prime+1)+(2p-2q^\prime-s_6^\prime)+(2s_6^\prime-1)=3p-q^\prime+s_6^\prime$ is used to color $n\leq 2p+3\times 2q^\prime +4\times (2p-2q^\prime-s_6^\prime) +8s_6^\prime=10p-2q^\prime +4s_6^\prime$ vertices. Note that in fact $k=3p-q^\prime+s_6^\prime$ only if there is a pair of segments of lengths $4$ and $5$ colored with three colors, as described above, but then the total number of vertices is at most $n\leq 10p-2q^\prime +4s_6^\prime - (8-4)-(8-5)=10p-2q^\prime +4s_6^\prime-7$. Otherwise $k=3p-q^\prime+s_6^\prime+1$. We have $3p-q^\prime+s_6^\prime\geq \frac{5}{2}p-\frac{1}{2}q^\prime+\frac{1}{2}+s_6^\prime$, so finally $k\geq \frac{n}{4}+1$ in either case.
\end{enumerate}

From Theorem \ref{thmPaths} it follows that for $n\geq 9$ we have $\chi_{md}(P_n)\leq \frac{n}{4}+1$, which concludes the proof.
\end{proof}

\begin{corollary}\label{corCycles}
Let $C=C_n$ be a cycle of order $n\geq 3$ with vertices $v_1, v_2, \ldots, v_{n-1}, v_{n}$. Then
$$
\chi_{md}(C)=
\begin{cases}
2,& \text{if }n\in\{4,6,8\},\\
3,& \text{if }n\in\{3,5,7,9,10\},\\
4,& \text{if }11\leq n\leq 13,\\
\left\lceil\frac{n}{6}\right\rceil+2,& \text{if }n\geq 14.
\end{cases}
$$
\end{corollary}

\begin{proof}

The fact that the given values are minimum possible follows immediately from the relation $\chi_{md}(C)\geq\chi(C)$ (for $n\leq 9$) and from Lemma \ref{lem_propertiesCycle} (for $n\geq 9$). 

To see that $\chi_{md}(C)$ does not exceed the given values one can observe that the proper colorings of $C_4$, $C_6$ and $C_8$ are also the majority dominator colorings, while for the other cycles one can use sequences of colors used for the paths of respective lengths defined in the proof of Theorem \ref{thmPaths}.
\end{proof}


\section{Final Remarks}\label{secFinal}

We introduced a new graph invariant, the majority dominator coloring number $\chi_{md}(G)$ and investigated some of its properties. Although at first glance the concept of majority dominator coloring is similar to the idea of dominator coloring, they are different enough to make it difficult to apply the same proof techniques. For that reasons, we needed to develop new ones.

We presented several general properties of the new parameter, in particular by connecting it with other graph invariants, like chromatic number $\chi(G)$, domination number $\gamma(G)$, independence number $\alpha(G)$ and matching number $\nu(G^\prime)$, where $G^\prime$ is a specific subgraph of $G$. This allowed us to derive the value of $\chi_{md}(G)$ for various classes of graphs and find the families of graphs satisfying chosen constraints imposed on $\chi_{md}(G)$. Obviously, there are many open problems to solve. Below we present those most interesting according to our opinion.

As it was observed in the beginning of Section \ref{secGeneral}, for every graph $G$ it holds $\chi_{md}(G)\leq \chi_d(G)$. Although one could expect that usually the inequality will be strict, in certain cases equality occurs (see e.g. the complete graphs $K_n$). For that reason our first question is as follows.

\begin{problem}
Characterize the graphs for which $\chi_{md}(G)=\chi_d(G)$.
\end{problem}

A similar problem is connected with Theorem \ref{m110}. As one can see even in the case of the products $C_n\circ K_1$ the bound is sometimes sharp and sometimes not (see Proposition \ref{thmCirc}). For that reason we would like to know the solution to the following problem.

\begin{problem}
Characterize the graphs for which $\chi_{md}(G)= \chi(G) +\left\lceil \frac{\alpha(G)}{2}\right\rceil-1$.
\end{problem}

Our last problem is connected with the product graphs $G\circ K_1$ for arbitrary $G$. In such case $\alpha(G\circ K_1)=n$, where $n$ is the order of $G$, and $\chi(G\circ K_1)=\chi(G)$ if $n\geq 2$. From Theorem \ref{m110} it follows that
$$
\chi_{md}(G\circ K_1)\leq \chi(G) +\left\lceil \frac{n}{2}\right\rceil-1
$$
On the other hand, by using a reasoning similar to the one from the proof of Proposition \ref{thmCirc}, we get
$$
\chi_{md}(G\circ K_1)\geq \max\left\{\chi(G),\left\lceil \frac{n}{2}\right\rceil+1\right\}.
$$
In particular, this allows us to make the following observation.
\begin{observation}
If $G$ is a bipartite graph of order $n$, then
$$
\chi_{md}(G\circ K_1)=\left\lceil \frac{n}{2}\right\rceil+1.
$$.
\end{observation}
However, the value of $\chi_{md}(G\circ K_1)$ remains unknown for general $G$.
\begin{problem}
Let $G$ be an arbitrary graph of order $n$. What is the value of $\chi_{md}(G\circ K_1)$? Can it be expressed in terms of $n$ and $\chi(G)$?
\end{problem}

\section*{Acknowledgment}
Marcin Anholcer was partially supported by the National Science Center of Poland under grant no. 2020/37/B/ST1/03298.

\end{document}